\documentclass[10pt]{article}
\usepackage{palatino}
\usepackage{amsmath}
\usepackage{amsthm}
\usepackage{amssymb}
\usepackage{mathdots} 
\usepackage{latexsym}
\usepackage[xcolor=pst]{pstricks}
\usepackage{graphicx, pst-plot, pst-node, pst-text, pst-tree}
\usepackage[small]{titlesec}
\usepackage[small,it]{caption}

\usepackage[dont-mess-around]{fnpct}

\usepackage{etoolbox}
\makeatletter
\patchcmd\maketitle{\setcounter{footnote}{0}}{}{}{}
\patchcmd\maketitle{%
  \renewcommand\thefootnote{\@fnsymbol\c@footnote}}{\AdaptNote\thanks\multthanks}{}{}
\patchcmd\maketitle{%
  \def\@makefnmark{\rlap{\@textsuperscript{\normalfont\@thefnmark}}}}{}{}{}
\makeatother

\usepackage{color}
\definecolor{lightgray}{rgb}{0.8, 0.8, 0.8}
\definecolor{darkgray}{rgb}{0.7, 0.7, 0.7}
\definecolor{darkblue}{rgb}{0, 0, .4}

\usepackage{hyperref}
\hypersetup{
        bookmarks=true,
        hyperfootnotes=false,
        colorlinks=true,
        linkcolor=darkblue,
        anchorcolor=darkblue,
        citecolor=darkblue,
        urlcolor=darkblue,
        pdfpagemode=UseThumbs,
        pdftitle={A simple proof of a theorem of Schmerl and Trotter for permutations},
        pdfsubject={Combinatorics},
        pdfauthor={Vincent Vatter},
}

\newcounter{todocounter}

\makeatletter
\newfont{\footsc}{cmcsc10 at 8truept}
\newfont{\footbf}{cmbx10 at 8truept}
\newfont{\footrm}{cmr10 at 10truept}
\renewenvironment{abstract}%
                {
                  \begin{list}{}%
                     {\setlength{\rightmargin}{1in}%
                      \setlength{\leftmargin}{1in}}%
                   \item[]\ignorespaces\begin{small}}%
                 {\end{small}\unskip\end{list}}

\newpagestyle{main}[\small]{
        \headrule
        \sethead[\usepage][][]
        {\sc Schmerl and Trotter for Permutations}{}{\usepage}}

\title{\sc A Simple Proof of a Theorem of Schmerl and Trotter for Permutations}

\newcommand{\footremember}[2]{%
	\footnote{#2}
	\newcounter{#1}
	\setcounter{#1}{\value{footnote}}%
}
\newcommand{\footrecall}[1]{%
	\footnotemark[\value{#1}]%
}

\author{%
\begin{tabular}{ccc}
Robert Brignall\footremember{epsrc}{Both authors were partially supported by the EPSRC Grant EP/J006130/1.}
&\rule{3pt}{0pt}&
Vincent Vatter\footrecall{epsrc}\textsuperscript{,}\footnote{Vatter was partially sponsored by the National Security Agency under Grant Number H98230-12-1-0207 and the National Science Foundation under Grant Number DMS-1301692.  The United States Government is authorized to reproduce and distribute reprints not-withstanding any copyright notation herein.}\\[-0.25ex]
\small Department of Mathematics and Statistics
&&
\small Department of Mathematics\\[-0.5ex]
\small The Open University
&&
\small University of Florida\\[-0.5ex]
\small Milton Keynes, England
&&
\small Gainesville, Florida USA\\[-1.5ex]
\end{tabular}
}

\titleformat{\section}
        {\large\sc}
        {\thesection.}{1em}{}

\date{}

\begin{document}
\maketitle

\pagestyle{main}

\begin{abstract}
When specialized to the context of permutations, Schmerl and Trotter's Theorem states that every simple permutation which is not a parallel alternation contains a simple permutation with one fewer entry. We give an elementary proof of this result.
\end{abstract}

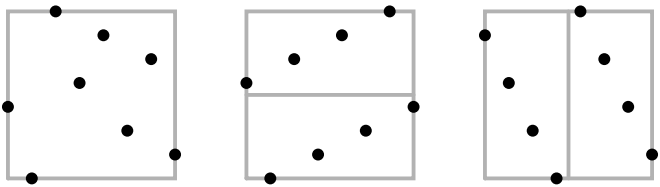
\begin{figure}[t]
\begin{center}
\begin{tabular}{ccccc}

\psset{xunit=0.0125in, yunit=0.0125in, runit=2.5\psxunit, linewidth=0.02in}
\begin{pspicture}(10,10)(80,80)
\psline[linecolor=darkgray]{c-c}(10,10)(80,10)(80,80)(10,80)(10,10)
\pscircle*(10,40){1}
\pscircle*(20,10){1}
\pscircle*(30,80){1}
\pscircle*(40,50){1}
\pscircle*(50,70){1}
\pscircle*(60,30){1}
\pscircle*(70,60){1}
\pscircle*(80,20){1}
\end{pspicture}

&\rule{3pt}{0pt}&

\psset{xunit=0.0125in, yunit=0.0125in, runit=2.5\psxunit, linewidth=0.02in}
\begin{pspicture}(10,10)(80,80)
\psline[linecolor=darkgray]{c-c}(10,10)(80,10)(80,80)(10,80)(10,10)
\psline[linecolor=darkgray]{c-c}(10,45)(80,45)
\pscircle*(10,50){1}
\pscircle*(30,60){1}
\pscircle*(50,70){1}
\pscircle*(70,80){1}
\pscircle*(20,10){1}
\pscircle*(40,20){1}
\pscircle*(60,30){1}
\pscircle*(80,40){1}
\end{pspicture}

&\rule{3pt}{0pt}&

\psset{xunit=0.0125in, yunit=0.0125in, runit=2.5\psxunit, linewidth=0.02in}
\begin{pspicture}(10,10)(80,80)
\psline[linecolor=darkgray]{c-c}(10,10)(80,10)(80,80)(10,80)(10,10)
\psline[linecolor=darkgray]{c-c}(45,10)(45,80)
\pscircle*(10,70){1}
\pscircle*(20,50){1}
\pscircle*(30,30){1}
\pscircle*(40,10){1}
\pscircle*(50,80){1}
\pscircle*(60,60){1}
\pscircle*(70,40){1}
\pscircle*(80,20){1}
\end{pspicture}

\end{tabular}
\end{center}
\caption{Three examples of simple permutations. The second and third are parallel alternations.}
\label{fig-example-simples}
\end{figure}

An \emph{interval} in the permutation $\pi$ (thought of in one-line notation) is a contiguous set of entries whose values also form a contiguous set. Every permutation of length $n$ has \emph{trivial} intervals of lengths $0$, $1$, and $n$, and permutations with only trivial intervals are called \emph{simple}. We take a graphical view of permutations, in which we identify a permutation $\pi$ with its \emph{plot}, the set of points  $(i,\pi(i))$ in the plane. Three examples of simple permutations are plotted in Figure~\ref{fig-example-simples}. Note that $1$, $12$, and $21$ are all simple and that there are no simple permutations of length three. A bit more examination shows that $2413$ and $3142$ are the only simple permutations of length four.

The second and third simple permutations in Figure~\ref{fig-example-simples} are called \emph{parallel alternations}. A parallel alternation is, formally, a permutation whose entries can be divided into two halves of equal length, either both increasing or both decreasing, such that the entries of the halves interleave perfectly. In particular, $12$, $21$, $2413$, and $3142$ are parallel alternations. While parallel alternations need not be simple ($1324$ is a nonsimple parallel alternation), from any parallel alternation we may obtain a simple permutation by removing at most two entries.

Specialized to permutations\footnote{The notions of intervals and simplicity extend naturally to all relational structures (though with different names, such as modules and primality). Schmerl and Trotter proved their result for simple, irreflexive, binary relational structures.}, the main result of Schmerl and Trotter is as follows.

\newtheorem*{schmerl-trotter-theorem}{The Schmerl-Trotter Theorem for Permutations~\cite{schmerl:critically-inde:}}
\begin{schmerl-trotter-theorem}
Every simple permutation which is not a parallel alternation contains an entry whose removal leaves a simple permutation.
\end{schmerl-trotter-theorem}

Schmerl and Trotter's theorem has found wide application both theoretically and practically. For example, it has been used in~\cite{albert:simple-permutat:,albert:inflations-of-g:} to show that certain classes of permutations are defined by a finite set of restrictions, and in Albert's \emph{PermLab} package~\cite{PermLab1.0} to efficiently generate the simple permutations in a permutation class.

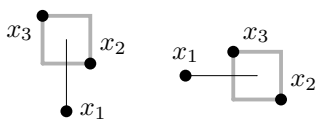
\begin{figure}[ht]
\centering
\begin{tabular}{ccc}

\psset{xunit=0.0125in, yunit=0.0125in, runit=2.5\psxunit, linewidth=0.005in}
\begin{pspicture}(0,0)(30,50)
\psline[linecolor=darkgray, linewidth=0.02in]{c-c}(10,30)(30,30)(30,50)(10,50)(10,30)
\psline[linecolor=black]{c-c}(20,10)(20,40)
\pscircle*(10,50){1}\uput[225](10,50){$x_3$}
\pscircle*(20,10){1}\uput[0](20,10){$x_1$}
\pscircle*(30,30){1}\uput[45](30,30){$x_2$}
\end{pspicture}

&\rule{3pt}{0pt}&

\psset{xunit=0.0125in, yunit=0.0125in, runit=2.5\psxunit, linewidth=0.005in}
\begin{pspicture}(0,-5)(50,30)
\psline[linecolor=darkgray, linewidth=0.02in]{c-c}(30,10)(50,10)(50,30)(30,30)(30,10)
\psline[linecolor=black]{c-c}(10,20)(40,20)
\pscircle*(10,20){1}\uput[90](10,20){$x_1$}
\pscircle*(30,30){1}\uput[45](30,30){$x_3$}
\pscircle*(50,10){1}\uput[45](50,10){$x_2$}
\end{pspicture}

\end{tabular}
\caption{In both figures, $x_1$ separates $x_2$ and $x_3$.}
\label{fig-separation}
\end{figure}

In this note we give a short, self-contained proof of the Schmerl-Trotter Theorem. We use only a few definitions in this proof. First, in order to simplify the discussion, we say that an entry of the simple permutation $\sigma$ is \emph{inessential} if its removal leaves a permutation which is still simple; otherwise, entries are \emph{essential}. The other concept we need is \emph{separation}: given a permutation $\sigma$ and entries $x_1$, $x_2$, and $x_3$, we say that $x_1$ \emph{separates $x_2$ and $x_3$} if $x_1$ lies between $x_2$ and $x_3$ either horizontally or vertically, but not both, as in Figure~\ref{fig-separation}. We extend this to sets of entries by saying that $x$ separates the entries $X$ if it lies outside the rectangular hull of $X$ and separates any two of its entries.

%
%
\begin{proof}
We prove the theorem by induction on the length of the simple permutation. It is vacuously true for permutations of length four, as both such simple permutations are parallel alternations, so suppose that $\sigma$ is a simple permutation of length at least five which is not a parallel alternation and that the theorem holds for all shorter permutations. As we are done otherwise, we assume throughout the proof that every entry of $\sigma$ is essential.

We begin by assuming, in order to eliminate this case, that removing any entry of $\sigma$ creates a minimal proper interval containing precisely two entries. Let $x_1$ be an arbitrary entry of $\sigma$ and suppose that $\sigma-x_1$ contains the interval $\{x_2,x_3\}$. It follows that $x_1$ must separate $x_2$ and $x_3$ in $\sigma$ because it is simple, and furthermore, $x_1$ must be the only entry which separates $x_2$ and $x_3$ in $\sigma$ because $\{x_2,x_3\}$ forms an interval when $x_1$ is removed. By applying one of the eight symmetries of the square, we may then assume that these entries are in the relative order of $231$, with $x_1$ on the left. This situation is depicted in Figure~\ref{fig-first-sep}. In this and all later figures, the gray areas indicate regions which cannot contain entries (because of the interval conditions).

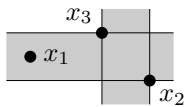
\begin{figure}[ht]
\centering
\begin{tabular}{c}
\psset{xunit=0.0125in, yunit=0.0125in, runit=2.5\psxunit, linewidth=0.005in}
\begin{pspicture}(0,0)(70,40)
\psframe[linecolor=lightgray,fillcolor=lightgray,fillstyle=solid](40,0)(60,40)
\psframe[linecolor=lightgray,fillcolor=lightgray,fillstyle=solid](0,10)(70,30)
\psline(0,10)(70,10)
\psline(0,30)(70,30)
\psline(40,0)(40,40)
\psline(60,0)(60,40)
\pscircle*(10,20){1}\uput[0](10,20){$x_1$}
\pscircle*(60,10){1}\uput[-45](60,10){$x_2$}
\pscircle*(40,30){1}\uput[135](40,30){$x_3$}
\end{pspicture}
\end{tabular}
\caption{The relative positions of $x_1$, $x_2$, and $x_3$.}
\label{fig-first-sep}
\end{figure}

Consider the doubleton interval in $\sigma-x_3$. Because $\sigma-x_3$ still contains $x_2$, this doubleton interval must contain at least one of $x_1$ or $x_2$, together with a new entry $x_4$. Figure~\ref{fig-doubleton-x4} shows the three possibilities: the doubleton interval can consist of $\{x_1,x_4\}$ with $x_4$ either to the left or right of $x_1$, or it can consist of $\{x_2,x_4\}$, but only if $x_4$ lies below and to the left of $x_2$.

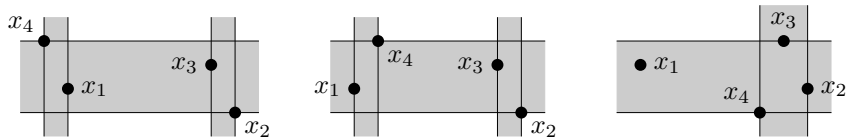
\begin{figure}[ht]
\centering
\begin{tabular}{ccccc}
\psset{xunit=0.0125in, yunit=0.0125in, runit=2.5\psxunit, linewidth=0.005in}
\begin{pspicture}(0,0)(100,50)
\psframe[linecolor=lightgray,fillcolor=lightgray,fillstyle=solid](10,0)(20,50)
\psframe[linecolor=lightgray,fillcolor=lightgray,fillstyle=solid](80,0)(90,50)
\psframe[linecolor=lightgray,fillcolor=lightgray,fillstyle=solid](0,10)(100,40)
\psline(0,10)(100,10)
\psline(0,40)(100,40)
\psline(10,0)(10,50)
\psline(20,0)(20,50)
\psline(80,0)(80,50)
\psline(90,0)(90,50)
\pscircle*(20,20){1}\uput[0](20,20){$x_1$}
\pscircle*(90,10){1}\uput[-45](90,10){$x_2$}
\pscircle*(80,30){1}\uput[180](80,30){$x_3$}
\pscircle*(10,40){1}\uput[135](10,40){$x_4$}
\end{pspicture}

&\rule{3pt}{0pt}&

\psset{xunit=0.0125in, yunit=0.0125in, runit=2.5\psxunit, linewidth=0.005in}
\begin{pspicture}(10,0)(100,50)
\psframe[linecolor=lightgray,fillcolor=lightgray,fillstyle=solid](20,0)(30,50)
\psframe[linecolor=lightgray,fillcolor=lightgray,fillstyle=solid](80,0)(90,50)
\psframe[linecolor=lightgray,fillcolor=lightgray,fillstyle=solid](10,10)(100,40)
\psline(10,10)(100,10)
\psline(10,40)(100,40)
\psline(30,0)(30,50)
\psline(20,0)(20,50)
\psline(80,0)(80,50)
\psline(90,0)(90,50)
\pscircle*(20,20){1}\uput[180](20,20){$x_1$}
\pscircle*(90,10){1}\uput[-45](90,10){$x_2$}
\pscircle*(80,30){1}\uput[180](80,30){$x_3$}
\pscircle*(30,40){1}\uput[-45](30,40){$x_4$}
\end{pspicture}

&\rule{3pt}{0pt}&

\psset{xunit=0.0125in, yunit=0.0125in, runit=2.5\psxunit, linewidth=0.005in}
\begin{pspicture}(0,0)(90,55)
\psframe[linecolor=lightgray,fillcolor=lightgray,fillstyle=solid](0,10)(90,40)
\psframe[linecolor=lightgray,fillcolor=lightgray,fillstyle=solid](60,0)(80,55)
\psline(0,10)(90,10)
\psline(0,40)(90,40)
\psline(60,0)(60,55)
\psline(80,0)(80,55)
\pscircle*(10,30){1}\uput[0](10,30){$x_1$}
\pscircle*(80,20){1}\uput[0](80,20){$x_2$}
\pscircle*(70,40){1}\uput[90](70,40){$x_3$}
\pscircle*(60,10){1}\uput[135](60,10){$x_4$}
\end{pspicture}
\end{tabular}
\caption{The three possible configurations after finding $x_4$.}
\label{fig-doubleton-x4}
\end{figure}

Next consider the doubleton interval in $\sigma-x_4$. In the leftmost case of Figure~\ref{fig-doubleton-x4}, there is only one possibility, as the doubleton could only be $\{x_3,x_5\}$ for a new entry $x_5$. Then we see that the doubleton interval in $\sigma-x_5$ must be $\{x_4,x_6\}$ for a new entry $x_6$. Continuing this process until we have run out of entries of $\sigma$ (as indicated on the left of Figure~\ref{fig-doubleton-x5}) yields the desired contradiction. The case where $x_4$ lies to the right of $x_1$ (the middle case in Figures~\ref{fig-doubleton-x4} and \ref{fig-doubleton-x5}) yields a similar contradiction. This leaves only the rightmost case of Figure~\ref{fig-doubleton-x4}. In this case we instead consider the doubleton interval in $\sigma-x_2$. As $x_1$ and $x_4$ must be separated by an entry outside the rectangular hull of $\{x_1,x_2,x_3,x_4\}$ in $\sigma$ (otherwise these four entries would form an interval in $\sigma$), this doubleton must consist of a new entry $x_5$ together with $x_3$, as shown on the right of Figure~\ref{fig-doubleton-x5}. However, in that case, $\sigma-x_4$ cannot contain a doubleton interval, which completes the contradiction to our assumption that $\sigma-x$ contains a doubleton interval for every entry $x$.

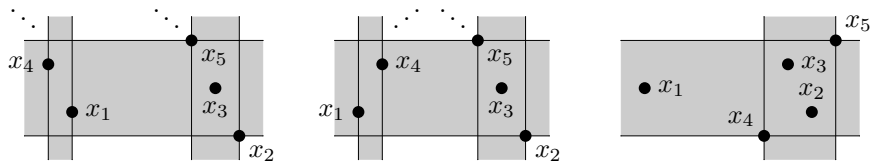
\begin{figure}[ht]
\centering
\begin{tabular}{ccccc}
\psset{xunit=0.0125in, yunit=0.0125in, runit=2.5\psxunit, linewidth=0.005in}
\begin{pspicture}(0,0)(100,60)
\psframe[linecolor=lightgray,fillcolor=lightgray,fillstyle=solid](10,0)(20,60)
\psframe[linecolor=lightgray,fillcolor=lightgray,fillstyle=solid](70,0)(90,60)
\psframe[linecolor=lightgray,fillcolor=lightgray,fillstyle=solid](0,10)(100,50)
\psline{c-c}(0,10)(100,10)
\psline{c-c}(0,50)(100,50)
\psline{c-c}(10,0)(10,60)
\psline{c-c}(20,0)(20,60)
\psline{c-c}(70,0)(70,60)
\psline{c-c}(90,0)(90,60)
\pscircle*(20,20){1}\uput[0](20,20){$x_1$}
\pscircle*(90,10){1}\uput[-45](90,10){$x_2$}
\pscircle*(80,30){1}\uput[-90](80,30){$x_3$}
\pscircle*(10,40){1}\uput[180](10,40){$x_4$}
\pscircle*(70,50){1}\uput[-45](70,50){$x_5$}
\uput[135](70,50){$\ddots$}
\uput[135](10,50){$\ddots$}
\end{pspicture}

&\rule{3pt}{0pt}&

\psset{xunit=0.0125in, yunit=0.0125in, runit=2.5\psxunit, linewidth=0.005in}
\begin{pspicture}(10,0)(100,60)
\psframe[linecolor=lightgray,fillcolor=lightgray,fillstyle=solid](20,0)(30,60)
\psframe[linecolor=lightgray,fillcolor=lightgray,fillstyle=solid](70,0)(90,60)
\psframe[linecolor=lightgray,fillcolor=lightgray,fillstyle=solid](10,10)(100,50)
\psline{c-c}(10,10)(100,10)
\psline{c-c}(10,50)(100,50)
\psline{c-c}(30,0)(30,60)
\psline{c-c}(20,0)(20,60)
\psline{c-c}(70,0)(70,60)
\psline{c-c}(90,0)(90,60)
\pscircle*(20,20){1}\uput[180](20,20){$x_1$}
\pscircle*(90,10){1}\uput[-45](90,10){$x_2$}
\pscircle*(80,30){1}\uput[-90](80,30){$x_3$}
\pscircle*(30,40){1}\uput[0](30,40){$x_4$}
\pscircle*(70,50){1}\uput[-45](70,50){$x_5$}
\uput[135](70,50){$\ddots$}
\uput[135](50,50){$\iddots$}
\end{pspicture}

&\rule{3pt}{0pt}&

\psset{xunit=0.0125in, yunit=0.0125in, runit=2.5\psxunit, linewidth=0.005in}
\begin{pspicture}(0,0)(100,60)
\psframe[linecolor=lightgray,fillcolor=lightgray,fillstyle=solid](0,10)(100,50)
\psframe[linecolor=lightgray,fillcolor=lightgray,fillstyle=solid](60,0)(90,60)
\psline{c-c}(0,10)(100,10)
\psline{c-c}(0,50)(100,50)
\psline{c-c}(60,0)(60,60)
\psline{c-c}(90,0)(90,60)
\pscircle*(10,30){1}\uput[0](10,30){$x_1$}
\pscircle*(80,20){1}\uput[90](80,20){$x_2$}
\pscircle*(70,40){1}\uput[0](70,40){$x_3$}
\pscircle*(60,10){1}\uput[135](60,10){$x_4$}
\pscircle*(90,50){1}\uput[45](90,50){$x_5$}
\end{pspicture}

\end{tabular}
\caption{The three possible configurations after finding $x_5$.}
\label{fig-doubleton-x5}
\end{figure}

The next case we consider is when there is an entry $x$ for which $\sigma-x$ contains a minimal proper interval, say $\iota$, of length four. By symmetry, we may further assume that $x$ lies to the left of $\iota$, which is itself a copy of $3142$, leaving us with three cases. These are depicted in Figure~\ref{fig-3142-cases}. In each of the three cases, the circled entry in the figure is a potential inessential entry of $\sigma$ which can only be essential if its removal creates a doubleton interval involving an entry of $\iota$ together with an entry $y$ which lies adjacent to $\iota$ both horizontally and vertically. Thus we may assume that the entry $y$ exists when handling these cases.

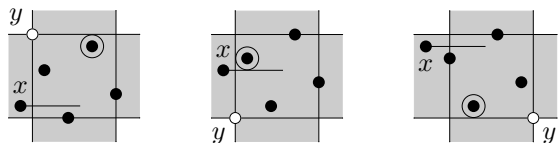
\begin{figure}[ht]
\centering
\begin{tabular}{ccccc}

\psset{xunit=0.0125in, yunit=0.0125in, runit=2.5\psxunit, linewidth=0.005in}
\begin{pspicture}(-5,0)(50,55)
\psframe[linecolor=lightgray,fillcolor=lightgray,fillstyle=solid](-5,10)(50,45)
\psframe[linecolor=lightgray,fillcolor=lightgray,fillstyle=solid](5,0)(40,55)
\psline{c-c}(-5,10)(50,10) 
\psline{c-c}(-5,45)(50,45) 
\psline{c-c}(5,0)(5,55) 
\psline{c-c}(40,0)(40,55) 
\pscircle*(10,30){1}
\pscircle*(20,10){1}
\pscircle*(30,40){1}
\pscircle*(40,20){1}
\psline(0,15)(25,15)
\pscircle*(0,15){1}\uput[90](0,15){$x$}
\pscircle[linewidth=0.005in](30,40){2}
\pscircle*[linecolor=white,linewidth=0.005in](5,45){1}
\pscircle[linewidth=0.005in](5,45){1}
\uput[135](5,45){$y$}
\end{pspicture}

&\rule{3pt}{0pt}&

\psset{xunit=0.0125in, yunit=0.0125in, runit=2.5\psxunit, linewidth=0.005in}
\begin{pspicture}(-5,-5)(50,50)
\psframe[linecolor=lightgray,fillcolor=lightgray,fillstyle=solid](-5,5)(50,40)
\psframe[linecolor=lightgray,fillcolor=lightgray,fillstyle=solid](5,-5)(40,50)
\psline{c-c}(-5,5)(50,5) 
\psline{c-c}(-5,40)(50,40) 
\psline{c-c}(5,-5)(5,50) 
\psline{c-c}(40,-5)(40,50) 
\pscircle*(10,30){1}
\pscircle*(20,10){1}
\pscircle*(30,40){1}
\pscircle*(40,20){1}
\psline(0,25)(25,25)
\pscircle*(0,25){1}\uput[90](0,25){$x$}
\pscircle[linewidth=0.005in](10,30){2}
\pscircle*[linecolor=white,linewidth=0.005in](5,5){1}
\pscircle[linewidth=0.005in](5,5){1}
\uput[225](5,5){$y$}
\end{pspicture}

&\rule{3pt}{0pt}&

\psset{xunit=0.0125in, yunit=0.0125in, runit=2.5\psxunit, linewidth=0.005in}
\begin{pspicture}(-5,-5)(55,50)
\psframe[linecolor=lightgray,fillcolor=lightgray,fillstyle=solid](-5,5)(55,40)
\psframe[linecolor=lightgray,fillcolor=lightgray,fillstyle=solid](10,-5)(45,50)
\psline{c-c}(-5,5)(55,5) 
\psline{c-c}(-5,40)(55,40) 
\psline{c-c}(10,-5)(10,50) 
\psline{c-c}(45,-5)(45,50) 
\pscircle*(10,30){1}
\pscircle*(20,10){1}
\pscircle*(30,40){1}
\pscircle*(40,20){1}
\psline(0,35)(25,35)
\pscircle*(0,35){1}\uput[270](0,35){$x$}
\pscircle[linewidth=0.005in](20,10){2}
\pscircle*[linecolor=white,linewidth=0.005in](45,5){1}
\pscircle[linewidth=0.005in](45,5){1}
\uput[-45](45,5){$y$}
\end{pspicture}

\end{tabular}
\caption{Up to symmetry, the three cases where $\iota$ has length four.}
\label{fig-3142-cases}
\end{figure}

In the first case, where $x$ separates the `$1$' and the `$2$' of $\iota$, the `$4$' of $\iota$ would be inessential unless its removal were to result in an interval involving the `$3$'. That possibility could only happen if there were an entry $y$ lying directly above and to the left of $\iota$ in $\sigma-x$. However in that case $y$ could not separate any set of entries not already separated by either $x$ or $\iota$ and thus $y$ would be inessential. In the second case there must be other entries of $\sigma$ (because $\sigma$ is not a parallel alternation), and so there must be an entry other than $y$ separating $\iota\cup\{x,y\}$. The only place such an entry can lie is vertically between $x$ and $y$, and this shows that $y$ is inessential. The third and final case is similar to the first.

Next suppose that $\sigma-x$ contains a minimal proper interval $\iota$ which is a parallel alternation of length $2\ell\ge 6$. By symmetry we may assume that $\iota$ is in the same relative order as $(\ell+1)1(\ell+2)2\cdots (2\ell)\ell$ and that $x$ lies either to the left of or above $\iota$. First suppose that $x$ lies above $\iota$. If $x$ lies horizontally between $\iota(1)$ and $\iota(2)$, then we have the situation shown on the left of Figure~\ref{fig-paralt-six}, and $\iota(2)$ is inessential. Otherwise, $x$ lies horizontally between two entries from the ``bottom half'' of $\iota$, say $\iota(2i)$ and $\iota(2i+2)$, and in this case the entry $\iota(2i+1)$ is inessential (the second picture in Figure~\ref{fig-paralt-six} shows an example of this).

\begin{figure}[ht]
\centering
\begin{tabular}{ccccccc}

\psset{xunit=0.0125in, yunit=0.0125in, runit=2.5\psxunit, linewidth=0.005in}
\begin{pspicture}(0,0)(90,95)
\psframe[linecolor=lightgray,fillcolor=lightgray,fillstyle=solid](0,10)(90,80)
\psframe[linecolor=lightgray,fillcolor=lightgray,fillstyle=solid](10,0)(80,95)
\psline{c-c}(0,10)(90,10) 
\psline{c-c}(0,80)(90,80) 
\psline{c-c}(10,0)(10,95) 
\psline{c-c}(80,0)(80,95) 
\pscircle*(10,50){1}
\pscircle*(20,10){1}
\pscircle*(30,60){1}
\pscircle*(40,20){1}
\rput{26.5650512}(50,70){$\dots$}
\rput{26.5650512}(60,30){$\dots$}
\pscircle*(70,80){1}
\pscircle*(80,40){1}
\psline(15,90)(15,45)
\pscircle*(15,90){1}
\uput[0](15,90){$x$}
\pscircle[linewidth=0.005in](20,10){2}
\end{pspicture}

&\rule{3pt}{0pt}&

\psset{xunit=0.0125in, yunit=0.0125in, runit=2.5\psxunit, linewidth=0.005in}
\begin{pspicture}(0,0)(90,95)
\psframe[linecolor=lightgray,fillcolor=lightgray,fillstyle=solid](0,10)(90,80)
\psframe[linecolor=lightgray,fillcolor=lightgray,fillstyle=solid](10,0)(80,95)
\psline{c-c}(0,10)(90,10) 
\psline{c-c}(0,80)(90,80) 
\psline{c-c}(10,0)(10,95) 
\psline{c-c}(80,0)(80,95) 
\pscircle*(10,50){1}
\pscircle*(20,10){1}
\pscircle*(30,60){1}
\pscircle*(40,20){1}
\rput{26.5650512}(50,70){$\dots$}
\rput{26.5650512}(60,30){$\dots$}
\pscircle*(70,80){1}
\pscircle*(80,40){1}
\psline(35,90)(35,45)
\pscircle*(35,90){1}
\uput[0](35,90){$x$}
\pscircle[linewidth=0.005in](30,60){2}
\end{pspicture}

&\rule{3pt}{0pt}&

\psset{xunit=0.0125in, yunit=0.0125in, runit=2.5\psxunit, linewidth=0.005in}
\begin{pspicture}(-10,-2.5)(90,90)
\psframe[linecolor=lightgray,fillcolor=lightgray,fillstyle=solid](-10,10)(90,80)
\psframe[linecolor=lightgray,fillcolor=lightgray,fillstyle=solid](10,0)(80,90)
\psline{c-c}(-10,10)(90,10) 
\psline{c-c}(-10,80)(90,80) 
\psline{c-c}(10,0)(10,90) 
\psline{c-c}(80,0)(80,90) 
\pscircle*(10,50){1}
\pscircle*(20,10){1}
\pscircle*(30,60){1}
\pscircle*(40,20){1}
\rput{26.5650512}(50,70){$\dots$}
\rput{26.5650512}(60,30){$\dots$}
\pscircle*(70,80){1}
\pscircle*(80,40){1}
\psline(0,15)(45,15)
\pscircle*(0,15){1}\uput[90](0,15){$x$}
\pscircle[linewidth=0.005in](30,60){2}
\end{pspicture}

&\rule{3pt}{0pt}&

\psset{xunit=0.0125in, yunit=0.0125in, runit=2.5\psxunit, linewidth=0.005in}
\begin{pspicture}(-5,-5)(90,90)
\psframe[linecolor=lightgray,fillcolor=lightgray,fillstyle=solid](-10,5)(90,80)
\psframe[linecolor=lightgray,fillcolor=lightgray,fillstyle=solid](5,-5)(80,90)
\psline{c-c}(-10,5)(90,5) 
\psline{c-c}(-10,80)(90,80) 
\psline{c-c}(5,-5)(5,90) 
\psline{c-c}(80,-5)(80,90) 
\pscircle*(10,50){1}
\pscircle*(20,10){1}
\pscircle*(30,60){1}
\pscircle*(40,20){1}
\rput{26.5650512}(50,70){$\dots$}
\rput{26.5650512}(60,30){$\dots$}
\pscircle*(70,80){1}
\pscircle*(80,40){1}
\psline(0,45)(45,45)
\pscircle*(0,45){1}\uput[120](0,45){$x$}
\pscircle[linewidth=0.005in](10,50){2}
\pscircle*[linecolor=white,linewidth=0.005in](5,5){1}
\pscircle[linewidth=0.005in](5,5){1}
\uput[225](5,5){$y$}
\end{pspicture}

\end{tabular}
\caption{Up to symmetry, the cases where $\iota$ is a parallel alternation of length at least six.}
\label{fig-paralt-six}
\end{figure}
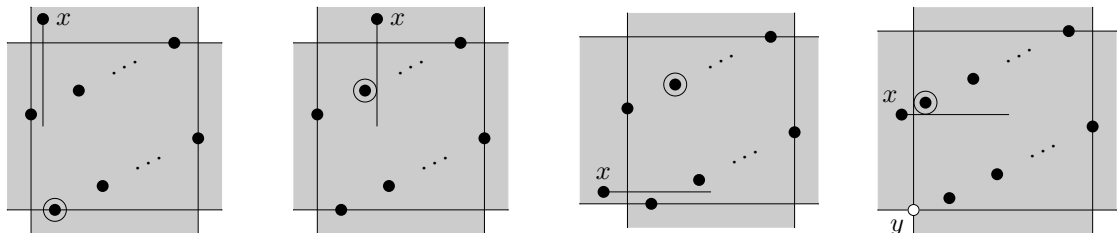

We may now suppose that $x$ lies to the left of $\iota$. As in the previous case, if $x$ lies vertically between two entries of the same ``half'' of $\iota$, say $\iota(i)$ and $\iota(i+2)$, then the entry $\iota(i+1)$ is inessential (the third picture in Figure~\ref{fig-paralt-six} shows an example of this). The only remaining case is when $x$ lies vertically between the first and last entries of $\iota$. In this case the only way $\iota(1)$ can be essential is if its removal creates an interval involving $\iota(2)$ and an entry $y$ lying immediately below and to the left of $\iota$ (as shown on the right of Figure~\ref{fig-paralt-six}). This case is similar to the analogous case where $\iota$ is of length four; $\sigma$ must have other entries because it is not a parallel alternation and in particular, it must contain an entry lying horizontally between $x$ and $y$ (because there is no other position in which $\iota\cup\{x,y\}$ can be separated). This implies that $y$ is inessential, thus completing the analysis of the case where the interval of $\sigma-x$ contains an interval which is a parallel alternation.

Having dealt with these special cases, we may now suppose that there is an entry $x_1$ such that $\sigma-x_1$ contains a minimal proper interval, say $\iota_1$, which is not a parallel alternation (and thus consists of at least five entries). We construct a sequence $x_1$, $x_2$, $\dots$ of essential entries and a sequence $\iota_1$, $\iota_2$, $\dots$ of intervals such that $\iota_k$ is a minimal proper interval of $\sigma-x_k$ for every $k$. Thus $x_k$ separates $\iota_k$ in $\sigma$, and the minimality of $\iota_k$ and simplicity of $\sigma$ imply that $\iota_k$ is itself simple. We are done by the above if any $\iota_k$ is a parallel alternation, so we may assume that each $\iota_k$ contains an inessential entry by induction, which we take to be $x_{k+1}$. 

\begin{figure}[ht]
\centering
\begin{tabular}{ccc}
\psset{xunit=0.0125in, yunit=0.0125in, runit=2.5\psxunit, linewidth=0.005in}
\begin{pspicture}(0,0)(70,75)
\psframe[linecolor=lightgray,fillcolor=lightgray,fillstyle=solid](0,10)(70,60)
\psframe[linecolor=lightgray,fillcolor=lightgray,fillstyle=solid](20,0)(60,75)
\psline{c-c}(0,60)(70,60) 
\psline{c-c}(0,10)(70,10) 
\psline{c-c}(20,0)(20,75) 
\psline{c-c}(60,0)(60,75) 
\pscircle*(20,20){1}
\pscircle*(40,30){1}
\pscircle*(50,10){1}
\pscircle*(60,50){1}
\psline(10,40)(40,40)
\pscircle*(10,40){1}
\uput[90](10,40){$x_k$}
\pscircle*(30,60){1}
\uput[30](30,60){$x_{k+1}$}
\rput[c](50,25){$\iota_k$}
\end{pspicture}

&\rule{3pt}{0pt}&

\psset{xunit=0.0125in, yunit=0.0125in, runit=2.5\psxunit, linewidth=0.005in}
\begin{pspicture}(0,0)(70,75)
\psframe[linecolor=lightgray,fillcolor=lightgray,fillstyle=solid](0,10)(70,60)
\psframe[linecolor=lightgray,fillcolor=lightgray,fillstyle=solid](20,0)(60,75)
\psline{c-c}(0,60)(70,60) 
\psline{c-c}(0,10)(70,10) 
\psline{c-c}(20,0)(20,75) 
\psline{c-c}(60,0)(60,75) 
\pscircle*(20,20){1}
\pscircle*(40,30){1}
\pscircle*(50,10){1}
\pscircle*(60,40){1}
\psline(10,50)(40,50)
\pscircle*(10,50){1}
\uput[-90](10,50){$x_k$}
\pscircle*(30,60){1}
\uput[30](30,60){$x_{k+1}$}
\rput[c](50,25){$\iota_k$}
\end{pspicture}
\end{tabular}
\caption{Examples of the two cases that can arise after choosing $x_{k+1}$.}
\label{fig-endgame}
\end{figure}
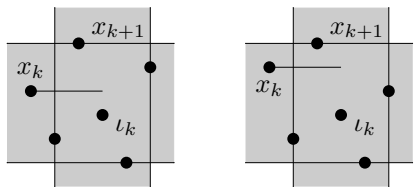

Next observe that because $\iota_k-x_{k+1}$ is simple and has length at least four, it does not have any entries on its corners. Therefore $\sigma-x_{k+1}$ cannot contain a minimal proper interval which includes entries from both inside and outside $\iota_k$. Now if $x_k$ separates $\iota_k-x_{k+1}$, as in the left of Figure~\ref{fig-endgame}, then $x_{k+1}$ is inessential for $\sigma$ and we are done. Otherwise, as shown on the right of Figure~\ref{fig-endgame}, $x_k$ must separate $x_{k+1}$ from the rest of $\iota_k$, i.e., from $\iota_k-x_{k+1}$. In this case $\iota_k-x_{k+1}$ is a proper interval of $\sigma-x_{k+1}$, so we take it to be $\iota_{k+1}$.

This process must terminate because $|\iota_{k+1}|=|\iota_k|-1$. Moreover, because each $\iota_k$ is simple and there are no simple permutations of length three, when this process does terminate it must be because we have either found an inessential entry of $\sigma$ or because some $\iota_k$ is a parallel alternation of length at least four, and in either case the theorem is proved.
\end{proof}

Now that we have proved that ``almost all'' simple permutations have an inessential entry, how many entries should we expect to be inessential? The answer is \emph{almost all of them}. This fact seems to have been first observed by Pierrot and Rossin~\cite{pierrot:simple-permutat:}, but we include a short sketch below.

It has been shown (see, for example, Corteel, Louchard, and Pemantle~\cite{corteel:common-interval:}) that the number of nontrivial intervals in a random permutation of length $n$ is asymptotically Possion distributed with mean $2$ (in fact this is true for nontrivial intervals of size $2$; the probability of a random permutation having a larger nontrivial interval tends to $0$ as $n\rightarrow\infty$). Therefore the number of simple permutations of length $n$ is asymptotic to $n!/e^2$. Now we double-count pairs $(\sigma,x)$ where $\sigma$ is a simple permutation of length $n+1$ and $x$ is an inessential entry of $\sigma$. On one hand, the number of such pairs is asymptotic to
$$
\frac{(n+1)!}{e^2}\cdot\mathbb{E}[\mbox{number of inessential entries}].
$$
On the other hand, the number of such pairs is equal to the number of pairs $(\sigma,\tau)$ where $\sigma$ and $\tau$ are simple, $\sigma$ has length $n+1$, and $\tau=\sigma-x$. Consider inserting a new entry $x$ into a simple permutation $\tau$ of length $n$. There naively $(n+1)^2$ different places to insert $x$. However, $2n$ of these places will create intervals of size two with entries of $\tau$ while the $4$ places on the corners will create an interval of size $n$. Each of the remaining $n^2-3$ places to insert $x$ yields a different permutation, and if inserted into one of those positions, $x$ cannot lie in an interval of size strictly between two and $n$ because $\tau$ has no proper intervals. Therefore the number of such pairs is asymptotic to
$$
\frac{n!}{e^2}\left(n^2-3\right),
$$
showing that the expected number of inessential entries is asymptotic to $n$, as desired.

\bigskip

\noindent{\bf Acknowledgements.}  We are grateful to Jay Pantone and the anonymous referees for their comments and corrections.

\bibliographystyle{acm}
\bibliography{../../refs}

\end{document}